\documentclass{amsart}
\usepackage[dvips]{color}
\usepackage{graphicx}

\usepackage{latexsym}
\usepackage{amssymb}
\usepackage{amsmath}
\usepackage{amsthm}
\usepackage{amscd}
\theoremstyle{plain}
\newtheorem{thm}{Theorem}
\newtheorem{lem}{Lemma}

\newtheorem{prop}{Proposition}

\theoremstyle{remark}
\newtheorem{rem}{Remark}

\newcommand{\Leb}{\ensuremath{\lambda}}

\newcommand{\LS}{\ensuremath{\underset{n=1}{\overset{\infty}{\cap}} \, {\underset{i=n}{\overset{\infty}{\cup}}}\,}}

\newtheorem{cor}{Corollary}
\author[J.\ Chaika]{Jon Chaika}\email{jonchaika@math.uchicago.edu}\address{Department of Mathematics, University of Chicago, 5734 S. University Avenue, Chicago, IL~60637, USA}
\begin{document}\title{Homogeneous Approximation for flows on translation surfaces}
\maketitle

Given a translation surface $M$ with a distinguished horizontal direction, we denote the flow in direction 2$\pi \theta$ to the horizontal by $F_{\theta}$. Let $d$ denote distance on $M$.

We say $M$ is a \emph{ Khinchin surface} if for every nonincreasing function $f: \mathbb{R}^+ \to \mathbb{R}^+$ such that $\int_C^{\infty}f(t) dt= \infty$  for any $C>0$ we have that for almost every $\theta\in [0,1)$ and $x\in M$, $d(F_{\theta}^t(x),x)<f(t)$ for arbitrarily large $t$.
\begin{thm}\label{main} All translation surfaces are Khinchin surfaces.
\end{thm}
\begin{cor}\label{main:cor} Let $\{a_i\}_{i=1}^{\infty}$ be a non-increasing sequence of positive real numbers with $\underset{i=1}{\overset{\infty}{\sum}}a_i=\infty$. For almost every IET $T$ and almost every point $x$ we have $\underset{n \to \infty}{\liminf}\, a_n^{-1}|T^nx-x| =0.$ That is, $|T^nx-x|<a_n$ for infinitely many $n$.
\end{cor}
We were motivated by: \begin{thm} \label{march} (Marchese \cite{LM}) Let $\{a_i\}_{i=1}^{\infty}$ be a decreasing sequence with divergent sum and with the additional property that $\{ia_i\}_{i=1}^{\infty}$ is decreasing. For almost every IET $T$ $$\delta \in \LS B(T^i(\delta'),a_i))$$ where $\delta$ and $\delta'$ are any discontinuities of $T$.
\end{thm}
\begin{rem}
Neither Corollary \ref{main:cor} or Theorem \ref{march} is stronger than the other. Corollary \ref{main:cor} allows more general sequences and treating almost every point. Marchese has recently adapted Theorem \ref{march} to treat almost every point. Additionally, Corollary \ref{main:cor} lets us treat a one parameter family obtained by fixing any translation surface and considering varying the direction of flow and examining first return to a transversal on a fixed translation surface. By treating each translation surface, it applies to billiards in rational polygons \cite{katzem}. It does not recover Theorem \ref{march}'s result that all discontinuities approximate each other.
\end{rem}

Let $M$ be a translation surface of genus $g$ and $s$ be the length of the shortest saddle connection on $M$. Let $\sigma=\frac 1 {2g-2}$. 
Let $S_{\sigma}(M)$ be the sequence of directions of periodic cylinders of volume at least $\sigma$ enumerated by increasing length. Let $\tilde{S}_{\sigma}(L)=$ $$\{(\theta,T): \theta \text{ is the direction of a periodic cylinder of length } T<L \text{ and volume }\sigma\}.$$


\begin{prop}\label{vorobets}(Vorobets)Let $m$ be the sum of the multiplicities of singularities. Let $c=2^{2^{4m}}\sqrt{s}$. For large enough $L$ we have $\underset{(\theta,T)\in \tilde{S}_{\sigma}(L)}{\cup}B(\theta,\frac{c}{TL})=S^1$.
\end{prop}
This is shown on \cite[page 14]{vor}.

The following proposition is a key result.
\begin{prop}\label{key} There exist $C_1,C_2,C_3>0$ such that for any interval $J$ and $N>C_3$ we have $\Leb(\underset{\theta \in \tilde{S}_{\sigma}(C_1N)\setminus \tilde{S}_{\sigma}(N)}{\cup} B(\theta, \frac 1 {(C_1N)^2})\cap J)>C_2\Leb(J)-\frac 2{(C_1N)}^2$.
\end{prop}
To prove this we establish a few other results.
\begin{lem}\label{flow} If $x$ is in a periodic cylinder with direction $\theta$,  length $L$ and area $\sigma$ then $x$ is not in a periodic cylinder of  direction $\phi$ and length $R$ for any $\phi \in B(\theta, \frac 1 { RL})$.
\end{lem}
\begin{proof}
Notice that, $d(F_{\theta}^{t\cos(\theta-\theta')}(x),F_{\theta'}^{t}(x))=|\sin(\theta-\theta')|t$ so long as no singularity lies on the geodesic connecting $F^s_{\theta}(x)$ to $F^s_{\theta'}(x)$ for $0\leq s\leq t$. The cylinder must have circumference $\frac{\sigma}{L}.$ If $x$ is in a periodic direction of length $R$ and direction $\phi$ then $F^s_{\phi}(x)=F^{s-R}_{\phi}(x)$ for all $s$. Choose $s$ so that $F^s_{\phi}(x)$ is just leaving the periodic cylinder in direction $\theta$ that $x$ is in.  So $F_{\phi}^{s-t}(x)$ is in this periodic cylinder for all $t<\frac{\sigma}{L}|\csc(\theta-\phi)|$ and in particular $F^{s-t}_{\phi}(x)\neq F_{\phi}^s(x)$.  So $R\geq \frac{\sigma}{L}|\csc(\theta-\phi)|$ which implies that $|\theta-\phi| \geq \frac{1}{RL}$.
\end{proof}
\begin{cor}\label{sum bound} For any interval $J$ we have $\Leb(\underset{(\theta,T)\in \tilde{S}_{\sigma}(L)}{\cup}B(\theta,\frac 1 {TL})\cap J) \leq \sigma^{-1}\Leb(J)$.
\end{cor}
\begin{proof} If $(\theta,T)\in \tilde{S}_{\sigma}(L)$ then there exists a periodic cylinder $U$ in direction $\theta$ with length $T$. If $x \in U$ then $x$ is not in another periodic cylinder of length less than $L$ and direction $\phi \in B(\theta , \frac 1 {TL}).$ So any direction in $J$ can be contained in at most $\sigma^{-1}$ of the $B(\theta',\frac 1 {LT})$ for $(\theta',L) \in \tilde{S}_{\sigma}(L)$.
\end{proof}
\begin{cor} There exists $C_5>0$ such that for all large enough $L$ we have $$|\{(\theta,T) \in \tilde{S}_{\sigma}(C_5L)\setminus \tilde{S}_{\sigma}(L):\theta \in J\}|>\frac{C_1} 2 \Leb(J)L^2.$$
\end{cor} 
\begin{proof} Proposition \ref{vorobets} and Corollary \ref{sum bound} imply that 
$$C_1\Leb(J)\leq\underset{(\theta,T ) \in \tilde{S}_{\sigma}(L), B(\theta,\frac 1{TL}) \cap J \neq \emptyset}{\sum}\frac{2}{TL}\leq  \sigma^{-1}(\Leb(J)+\frac{2}{L}).$$
Therefore for large enough $L$ we have $\underset{(\theta,T ) \in \tilde{S}_{\sigma}(L)}{\sum}\frac{2}{T(3C_1^{-1}\sigma^{-1} )L}<\frac{C_1\Leb(J)}2$. So 
$$\underset{(\theta,T ) \in \tilde{S}_{\sigma}(3C_1^{-1}\sigma^{-1} L)\setminus \tilde{S}_{\sigma}(L), B(\theta,\frac 1{TL}) \cap J \neq \emptyset} {\sum}\frac{2}{3C_1^{-1}\sigma^{-1}TL}\geq \frac{C_1\Leb(J)}2.$$ 
The corollary follows.
\end{proof}
\begin{proof}[Proof of Proposition \ref{key}] This is a consequence of Lemma \ref{flow} and the previous corollary because for fixed $k>1$ and  $(\theta,T) \in \tilde{S}_{\sigma}(kR)\setminus \tilde{S}_{\sigma}(R)$ we have $TkR$ is proportional to $(kR)^{2}$.
\end{proof}
We now state a useful theorem which follows from \cite{ber russian} and is phrased as in \cite[Theorem 1]{BC}.
\begin{thm}\label{bc} If $\{x_i\}_{i=1}^{\infty} \subset [0,1]$ is a sequence with the property that there exists $c>0$ such that for any interval $J$ we have 
$$\underset{N \to \infty}{\liminf} \Leb(\underset{i=1}{\overset{N}{\cup}}B(x_i,s_i)\cap J)>c\Leb(J)$$ then for any non increasing sequence of positive real numbers $r_1 \geq r_2\geq ...$ with $\sum a_i=\infty$ we have $\Leb(\LS B(s_i, r_i))=1$.
\end{thm}
 
Observe that Proposition \ref{key} shows that $S_{\sigma}(M)$ satisfies the hypothesis of this theorem. So by fact that the number of cylinders is quadratic in the length (that is, the $k^{th}$ element has length proportional to $\sqrt{k}$) we obtain the following corollary.
\begin{cor}\label{use} If If $f$ is a non-increasing positive function such that $\int _C ^{\infty} f(\sqrt{t})dt=\infty$ for any $C>0$ then for almost every $\phi$ and any $L_0>0$ there exists $(\theta,T) \in S_{\sigma}(\infty)$ with $T>L_0$ and $|\phi-\theta|<f(T)$. 
\end{cor}
\begin{lem} \label{translation} Let $\theta$ be the direction of a cylinder of length $L$ and volume $a$. If 
$|\theta-\phi|\leq  \frac{ \epsilon}L$ then for a set of points of measure at least $\frac {a}{4}$ we have $d(F_{\phi}^{L|\sec(\theta-\phi)|}x,x)<2\epsilon$.
\end{lem}
\begin{proof} $d(F_{\theta}^{t\cos(\theta-\theta')}(x),F_{\theta'}^{t}(x)) =|\sin(\theta-\theta')|t$ so long as no singularity lies on the geodesic connecting $F^s_{\theta}(x)$ to $F^s_{\theta'}(x)$ for $0\leq s\leq t$. Therefore $d(x,F_{\phi}^{L\sec(\theta-\theta')})=d(F_{\theta}^L(x),F_{\phi}^{L\sec(\theta-\theta')}(x))=L\sec(\theta-\theta')\sin(\theta-\phi)$ for any $x$ such that $F^s_{\phi}(x)$ is contained in the cylinder. The lemma follows by the fact that $|\sin(\theta')|\leq |\theta'|$.
\end{proof}
\begin{prop}\label{close} If $f$ is a non-increasing positive function $\int_C^{\infty}tf(t)dt=+\infty$ for all $C>0$ then for almost every $\phi$ and any $L_0$ we have $|\theta-\phi|<f(L)$ where $\phi$ is the direction of a periodic cylinder of volume at least $\sigma$ and length $L>L_0$.
\end{prop}
\begin{proof}If $f$ is a non-increasing positive function then $\int _C ^{\infty} f(\sqrt{t})dt=\infty$ iff $\int _C ^{\infty} tf(t)dt=\infty$.
 It suffices to prove the discrete version. Observe that $\underset{i=1}{\overset{\infty}{\sum}}ia_i=\infty$ if and only if
$\underset{i=1}{\overset{\infty}{\sum}}a_{\lfloor \sqrt{i} \rfloor}=\infty$.
Indeed, let $M\geq 2$ and notice
$
\underset{i=1}{\overset{\infty}{\sum}} M^{2(i-1)}a_{M^i}\leq \underset{i=1}{\overset{\infty}{\sum}}ia_i\leq \underset{i=0}{\overset{\infty}{\sum}} M^{2(i+1)}a_{M^{i}},
$
and, on the other hand,
$\underset{i=1}{\overset{\infty}{\sum}} M^{2(i-1)}a_{M^i}\leq \underset{i=1}{\overset{\infty}{\sum}}a_{\lfloor \sqrt{i} \rfloor}\leq \underset{i=0}{\overset{\infty}{\sum}} M^{2(i+1)}a_{M^{i}}.
$ The proposition now follows by Corollary \ref{use}.
\end{proof}
\begin{proof}[Proof of Theorem \ref{main}] Because $f$ is non-increasing, if the flow in direction $\phi$ is ergodic it suffices to show that $\underset{t \to \infty}{\liminf}f(t)^{-1}d(F_{\phi}^t(x),x)=0$ for a positive measure set of $x$. $F_{\phi}$ is ergodic for Lebesgue almost every $\phi$ by \cite{KMS}. By Lemma \ref{translation} it suffices to show that for any $\epsilon>0$ there are arbitrarily large $L_0$ such that there is a  periodic cylinder of length $L>L_0$ and at least a fixed positive area in direction $\theta$ where $|\phi-\theta|<\frac {\epsilon f(L)}{L}$. This is true by Proposition \ref{close}. 
\end{proof}

\section{Acknowledgments}
I would like to thank M. Boshernitzan and L. Marchese.

\end{document}